\documentclass[reqno,12pt]{amsart} 
\usepackage{mathrsfs}
\usepackage{amssymb,amsmath,amsthm,amsfonts}
\usepackage[dvipsnames]{xcolor}
\usepackage[
		colorlinks=true,
		linkcolor=blue,
		citecolor=ForestGreen,
		urlcolor=black
		]
{hyperref}
\usepackage[inline]{enumitem}  
\usepackage[margin=1in]{geometry}
\let\oldsquare\square
\usepackage{mathabx}
\usepackage{graphicx}
\usepackage{cases}
\usepackage{cite}

\usepackage{tikz}\usepackage{subfigure}
\usetikzlibrary{decorations.pathreplacing,calligraphy}
\newcommand{\n}[1]{\left\Vert#1\right\Vert}
\newcommand{\abs}[1]{\left\vert#1\right\vert}
\newcommand{\set}[1]{\left\{#1\right\}}

\newcommand{\R}{\mathbb R}
\newcommand{\rd}{\R^d}
\newcommand{\C}{\mathbb C}
\newcommand{\Z}{\mathbb Z}

\newcommand{\N}{\mathbb N}

\newcommand{\F}{ {\mathcal F} }
\newcommand{\h}{ {\mathcal H} }

\newcommand{\1}{\mathbf{1}}
\newcommand{\m}{\mathcal{M}}
\newcommand{\mh}{\widehat{\m}}
\newcommand{\td}{\mathbb{T}^d}
\newcommand{\zd}{\Z^d}

\newcommand{\eps}{\varepsilon}

\newcommand{\s}{2\sigma+1}

\newcommand{\Conv}{\mathop{\scalebox{1.5}{\raisebox{-0.2ex}{$\ast$}}}}

\def\TagOnRight

\theoremstyle{plain}
\newtheorem{theorem}{Theorem} [section]
\newtheorem{lemma}[theorem]{Lemma}

\theoremstyle{remark}
\newtheorem{remark}[theorem]{Remark}

\theoremstyle{definition}
\newtheorem{definition}[theorem]{Definition}

\newtheorem*{A1}{{\bf Theorem A}}
\newtheorem*{A2}{{\bf Theorem B}}

\def\({\left(}
\def\){\right)}
\def\<{\left\langle}
\def\>{\right\rangle}

\numberwithin{equation}{section}

\begin{document}
\title[Norm inflation  for NLW ]
{Norm inflation with infinite loss of regularity at general initial data for nonlinear wave equations  in Wiener amalgam and Fourier amalgam spaces}  
\author[D. Bhimani]{Divyang G. Bhimani}
\address{Department of Mathematics, Indian Institute of Science Education and Research, Dr. Homi Bhabha Road, Pune 411008, India}
\email{divyang.bhimani@iiserpune.ac.in}
\author[S. Haque] {Saikatul Haque}
\address{Harish-Chandra Research Institute, Chhatnag Road, Jhunsi, 
Prayagraj (Allahabad) 211019, India}
\email{saikatulhaque@hri.res.in}
\thanks{} 

\maketitle
\begin{abstract}
  We study the strong  ill-posedness (norm inflation with infinite loss of regularity) for the nonlinear  wave equation at every initial data in Wiener amalgam and  Fourier amalgam spaces with negative regularity. In particular these spaces contain Fourier-Lebesgue, Sobolev and some modulation spaces. The equations are posed on $\mathbb R^d$ and on torus $\td$ and  involve a smooth power nonlinearity. Our results are sharp with respect to  well-posedness results of B\'enyi and Okoudjou (2009) and  Cordero and Nicola (2009) in the Wiener amalgam and modulation space cases.  In particular, we also complement  norm inflation result of Christ, Colliander and Tao  (2003)  and Forlano and  Okamoto (2020)  by establishing  infinite loss of regularity in the aforesaid spaces.    
\end{abstract}
\section{Introduction}
We study strong ill-posedness for  nonlinear wave (NLW) equations of the form
\begin{equation} \label{nlw}
    \begin{cases}
    \partial_t^2 u -  \Delta u = \pm  {u^{\rho}(\overline{u})^{\sigma-\rho}}\\
    ( u(0, \cdot),u_t(0,\cdot))=(u_0,u_1)
    \end{cases}
    \end{equation}
 $(t, x)\in \mathbb R \times \mathbb \m$,  where $\m=\rd$ or $\td$, 
{$\sigma \in \mathbb N$, $\rho\in\N\cup\{0\}$ with $\sigma\geq\max(\rho,2)$},  in Fourier amalgam and Wiener amalgam spaces.   
We recall:
\begin{definition}[Hadamard's well/ill posedness]\label{wpd} \
\begin{itemize}
\item The Cauchy problem \eqref{nlw} is called locally well-posed from $\mathcal{X}=X_0\times X_1$ to $Y$ if for every bounded set $B\subset\mathcal{X}$, there exist $T>0$ and a Banach space  $X_{T} \hookrightarrow C([0,T], Y)$ such that
(i) for all $\vec{u}_0:=(u_0,u_1) \in B$,  \eqref{nlw}  has a unique solution $u\in X_T$ with  $(u(0,\cdot),\partial_t u(0,\cdot))=\vec{u}_0$
(ii)  the solution map $\vec{u}_0 \mapsto u$ is continuous from $(B, \|\cdot\|_{\mathcal{X}})$ to $C([0,T],Y).$ 
In particular, when $\mathcal{X}=Y_s\times Y_{s-1}$ (in our context, $Y_s=\widehat{w}_s^{p,q}$ or $W_s^{2,q}$ to be defined below) and $Y=Y_s$, we say \eqref{nlw} is  locally well-posed in $Y_s$.
\item  The  Cauchy problem \eqref{nlw} is \textit{ill-posed} if the solution map is not continuous.
\item We say \textit{norm inflation} (NI) occurs at $\vec{u}_0\in \mathcal{X} $ in $\mathcal{X}$ (in short $NI_{\mathcal{X}}(\vec{u}_0)$) if 
given $\epsilon>0$, there exist $\vec{v}_0\in\mathcal{X},0<t<\epsilon$ with $\n{\vec{u}_0-\vec{v}_0}_{\mathcal{X}}<\epsilon$ such that for the solution $v$ to \eqref{nlw} corresponding to the data $\vec{v}_0$ one has $\n{v(t)}_{X_0}>\epsilon^{-1}$.   If $\mathcal{X}=\mathcal{X}_s=Y_s\times Y_{s-1}$, we write $NI_s(\vec{u}_0)$ for $NI_{\mathcal{X}}(\vec{u}_0)$ if there is no confusion.
\item If $NI_{\mathcal{X}}(0)$ occurs, we say  (mere) norm inflation (at zero) occurs to \eqref{nlw} in $\mathcal{X}$.
\item We say NI occurs with \textit{infinite loss of regularity} at ${\vec{u}_0}$ in $\mathcal{X}_s$ if for any given $\theta\in\R,\epsilon>0$,  there exist $\vec{v}_0\in\mathcal{X}_s,0<t<\epsilon$ with $\n{\vec{u}_0-\vec{v}_0}_{\mathcal{X}_s}<\epsilon$ such that for the  solution $v$ to \eqref{nlw} corresponding to the data $\vec{v}_0$ one has $\n{{v}(t)}_{X_\theta}>\epsilon^{-1}$.
\item If the conditions in the preceding definition occurs for $\theta$ on a proper subset of $\R$, we say NI occurs with \textit{finite loss of regularity} at ${\vec{u}_0}$ in $\mathcal{X}_s$.
\end{itemize}
\end{definition}
Recently in \cite{JH},   Oh and Forlano have introduced Fourier amalgam spaces $\widehat{w}^{p,q}_s(\m)$  (with $1\leq p, q \leq \infty,  s\in \mathbb R)$:
$$\widehat{w}^{p,q}_s(\m)=\left \{ f\in \mathcal{S}'(\m): \|f\|_{\widehat{w}^{p,q}_s}=   \left\| \left\lVert \chi_{n+Q_1}(\xi)\F f(\xi)\right\rVert_{L_\xi^p(\mh)}  \langle n \rangle^s \right\|_{\ell^q_n(\zd)}< \infty \right\},
$$
where $Q_1=(-\frac{1}{2},\frac{1}{2}]^d,$  $\mathcal{F}$ denotes the Fourier transform and  $\langle \cdot \rangle^s= (1+|\cdot|^2)^{s/2}.$  Here $\widehat{\m}$ denotes the Pontryagin dual of $\m,$ i.e.  $\widehat{\m}= \rd$ if $\m=\rd$ and $\widehat{\m}=\mathbb Z^d$ if $\m=\td$.  See \cite[Part II]{RuzhanskyTurunen},  among others, for  
 details on space $\mathcal{S}'(\m).$  The Fourier-Lebesgue spaces $\mathcal{F}L^q_s(\m)$ is defined by 
$$\mathcal{F}L^q_s(\m)=\left\{f\in \mathcal{S}'(\m): \mathcal{F}f \  \langle  \cdot  \rangle^s  \in L^q (\widehat{\m})\right \}.$$ 
 On the  other hand,  the modulation $M^{p,q}_s(\m)$ and Wiener amalgam $W^{p,q}_s(\m)$ 
 spaces were introduced by Feichtinger in early  1980's  in \cite{feichtinger1983modulation}.
To recall their definitions, let   $\rho \in \mathcal{S}(\R^d),$  $\rho: \R^d \to [0,1]$  be  a smooth function satisfying   $\rho(\xi)= 1 \  \text{if} \ \ |\xi|_{\infty}\leq \frac{1}{2} $ and $\rho(\xi)=
0 \  \text{if} \ \ |\xi|_{\infty}\geq  1.$ Let  $\rho_n$ be a translation of $\rho,$ that is,
$ \rho_n(\xi)= \rho(\xi -n), n \in \Z^d$
and denote 
$\sigma_{n}(\xi)=
  \frac{\rho_{n}(\xi)}{\sum_{\ell\in\Z^{d}}\rho_{\ell}(\xi)},  n
  \in \Z^d.$
Then the frequency-uniform decomposition operators can be defined by 
\[\oldsquare_n = \mathcal{F}^{-1} \sigma_n \mathcal{F}. \]
Now 
the modulation $M^{p,q}_s(\m)$ and  and Wiener amalgam spaces $W^{p,q}_s(\m)$, (with $1\leq p, q \leq \infty, s\in \mathbb R$) are defined by the norms:  
\begin{equation*}
\|f\|_{M^{p,q}_s(\m)}=   \left\| \left\lVert
  \oldsquare_nf\right\rVert_{L_x^p(\m)} \langle n \rangle ^{s} \right\|_{\ell^q_n(\zd)}  \quad \text{and} \quad \|f\|_{W^{p,q}_s(\m)}=   \left\| \left\lVert
  \oldsquare_nf \cdot \langle n \rangle ^{s}  \right\|_{\ell^q_n(\zd)}  \right\rVert_{L_x^p(\m)} .
\end{equation*}
 See Remark \ref{stft} for equivalent characterization via short-time Fourier transform (STFT) of these spaces. It is known that 
\begin{equation*} 
   \widehat{w}_s^{p,q}(\m)= \begin{cases}
   \mathcal{F}L_s^q(\m) \quad \textit{if}  \quad p=q\\
   M^{2,q}_s(\m)\quad \textit{if}  \quad p=2\\
   H^{s}(\m)=M^{2,2}_s(\m)=W^{2,2}_s(\m) \quad \textit{if} \quad p=q=2\\
   \mathcal{F}L_s^q(\m)=M_s^{p,q}(\m)=W^{p,q}_{s}(\m) \quad \textit{if}    \  \m =\td.
    \end{cases}
    \end{equation*}
(See \cite[Section 5]{CV1} and \cite[Proposition 11.3.1]{grochenig2013foundations} for details.)
In last two decades, modulation $M^{p,q}_s$ and Wiener amalgam $W^{p,q}_s$ spaces have been extensively studied in  PDEs, see  e.g. \cite{benyi2007unimodular, ruzhansky2012modulation, wang2011harmonic, wang2007global, benyi2009local, baoxiang2006isometric, bhimani2020hartree, bhimani2016functions, bhimani2020norm, bhimani2021strong, Elenawave, Reich}.  
  In particular,  we have  at least good  local well-posedness theory for NLW  \eqref{nlw} in these spaces.  We summarize them in the following:    
\begin{A1}[well-posedness]\label{hwp} Let  $1\leq p, q \leq \infty,$ $q'>(\sigma-1)d$ where  $q'$ is the H\"older conjugate of $q$, and let $Y_s=W^{p,q}_s(\rd)$ or $M^{p,q}_s(\mathbb R^d)$ or $\F L^1(\rd)$.
\begin{enumerate}
\item \cite{benyi2009local,  kassob, bhimani2016functions, Elenawave, Reich,forlano2019remark}   NLW \eqref{nlw} is locally well-posed in $Y_s$ for $ s\geq 0.$
\item \cite{TaoW, Lindblad}   NLW \eqref{nlw} is locally well-posed in $H^s(\rd)$ for $s\geq s(\sigma,d), d\geq2,$  some $s(\sigma,d)>0$. 
\end{enumerate}
\end{A1}
To the best of authors knowledge there is no local well-posedness results for \eqref{nlw} for in $\widehat{w}_s^{p,q}$ spaces (except the cases $p= 2$ or $p=q=1,s=0$).   
On the other hand,  there is extensive literature on ill-posedness for \eqref{nlw}.
We summarize some of them which are most suitable in our context.
\begin{A2}[ill-posedness] \
\begin{enumerate}
\item  \cite[Theorem 6]{christ2003ill}   Norm inflation  at zero occurs  for  \eqref{nlw} with nonlinearity $\pm\abs{u}^2u$ 
 in $H^{s}(\mathbb R)$ for $-1/2<s<0.$

\item  \cite[Theorem 1.7]{forlano2019remark} NI at general initial data for   \eqref{nlw} with  nonlinearity $\pm u^{\sigma}$ (the case $\rho=\sigma$ in our setting) occurs in $\widehat{w}^{p,q}_s(\rd)$   $(1\leq p,q \leq \infty)$ for $s<0.$
\item   \cite{Lebeau}  NI  with  finite loss of regularity for  \eqref{nlw} with nonlinearity $-u^{\sigma} (\sigma>3)$   in $H^{s}(\mathbb R^d)$ with some positive $s$ (precisely for $1<s< d/2- 2/(\sigma-1), d\geq 3).$
\end{enumerate}
\end{A2}
The aim  of this note is to  complement  positive results (Theorem A)  by  establishing strong ill-posedness for \eqref{nlw} in these spaces with negative regularity $s<0.$  We also complement  ill-posedness results (Theorem B) by  exhibiting infinite loss of regularity at general initial data.  We now state our main result.   
\begin{theorem}\label{dst} Assume that $1\leq p, q\leq \infty, s<0$ and let
\[ X^{p,q}_{s}(\m)=  \begin{cases}\widehat{w}_s^{p,q}(\rd) \text{ or }    W_s^{2,q}(\rd)\text{ for }\m=\rd\\
\F L_s^q(\td)      \text{ for }\m=\td 
\end{cases} \text{and} \quad \mathcal{X}_s^{p,q}(\m)=X_s^{p,q}(\m)\times X_{s-1}^{p,q}(\m).
  \] 
Then norm inflation with infinite loss of regularity  occurs to \eqref{nlw} at each element in $\mathcal{X}_s^{p,q}(\m)$:
  For any $\theta\in\R$, $\varepsilon>0$ and $\vec{u}_0\in\mathcal{X}_s^{p,q}(\m)$ there exist $\vec{u}_{0,\varepsilon} \in \mathcal X _s^{p,q}(\m)$ and  $T>0$ satisfying  
 \[ \|\vec{u}_0-\vec{u}_{0,\varepsilon}\|_{\mathcal{X}_s^{p,q}}< \varepsilon, \quad 0<T< \varepsilon \] 
 such that the corresponding smooth solution $u_\varepsilon$ to $\eqref{nlw}$ with data $\vec{u}_{0,\varepsilon}$ exists on $[0,T]$ and 
 \[ \|u_\varepsilon(T)\|_{X_\theta^{p,q}}> \frac{1}{\varepsilon}
 .\] 
 In particular, for any $T>0,$ the solution map $X_s^{p,q} \times X_{s-1}^{p,q}  \ni (u_0, u_1)\mapsto u \in C([0, T], X_\theta^{p,q})$ for  \eqref{nlw}  is discontinuous everywhere in $X_s^{p,q} \times X_{s-1}^{p,q}$ for all $\theta\in\mathbb R$.
\end{theorem}

For $q'>(\sigma-1) d$, Theorem \ref{dst} is sharp in the sense that  \eqref{nlw} is strongly ill-posed in $W^{2,q}_s(\rd)$, $M^{2,q}_s(\rd)$ for $s<0$ while by Theorem A it is locally well-posed for $s\geq 0.$  In fact,  even \textit{mere} ill-posedness in $W^{2,q}_s(\rd)$ is completely new as it is not cover by Theorem B. The particular case $X^{p,q}_{s}=\widehat{w}_s^{p,q}$ of Theorem \ref{dst} recover Theorem B (2) and further reveal that even worse situation occurs as  the infinite loss of regularity (everywhere in $\widehat{w}_s^{p,q}\times \widehat{w}_{s-1}^{p,q}$)  is present in all dimensions.  Theorem \ref{dst}
 also complements Theorem B (3) by taking $\rho=\sigma$ in \eqref{nlw} 
 by establishing infinite loss of regularity.

The circle of ideas  to  establish ill-posedness, via Fourier analytic approach, is originated  from the abstract argument \cite[Theorem 2 and Section 3]{BejenaruTao} of  Bejenaru and Tao  in the context of  quadratic nonlinear Schr\"odinger equation (NLS).  In fact,  the idea is to  rewrite the solution of \eqref{nlw} as a  power series expansion (Lemma \ref{uc}) in terms of Picard iterations.  It is then sufficient to establish the discontinuity for one Picard iterate to get the discontinuity of solution map at zero.   This method is further developed by  Iwabuchi and Ogawa  \cite{IwabuchiOgawa} to establish stronger phenomena of NI for NLS.  Later Kishimoto  \cite{kishimoto2019remark}   establish NI for NLS with more general  nonlinearity\footnote{Precisely this nonlinarity:  $\sum_{j=1}^n \nu_ju^{\rho_j}(\bar{u})^{\sigma_j-\rho_j}$ where $\nu_j\in\C$, $\sigma_j \in \mathbb N$, $\rho_j\in\N\cup\{0\}$ with $\sigma_j\geq\max(\rho_j,2)$.}  on the special domain  $\mathbb R^{d_1} \times \mathbb T^{d_2}, d=d_1+d_2.$   The idea is to show that one term in the series exhibits instability and dominates all the other terms (Lemmata \ref{d2'}, \ref{dt2'} and \ref{d1}, \ref{dt1'}) after adding a perturbation $\vec{\phi}_{0,N}$ to the data $\vec{u}_0$, see \eqref{data}.   We shall  notice that the existence time $T>0$ is allowed to shrink for the purpose of establishing norm inflation while in \cite{BejenaruTao}  it is fixed and uniform with respect to the initial data. See \cite[Section 2]{kishimoto2019remark}  for detail discussion on the approach. 
 On the other hand, Oh  \cite{OhT} use power series expansion indexed by trees to establish NI at general initial data for cubic NLS  on $\m$.
  Forlano and Okomoto \cite{forlano2019remark} use this Fourier analytic approach, following presentation from \cite{OhT},   for NLW \eqref{nlw} (in the particular case $\rho=\sigma$) to establish norm inflation at general initial data.
  
  In  addition to these ideas,  in the present paper,  we fix the size of the support of  perturbation $\vec{\phi}_{0,N}$ in \eqref{data} of  the initial data on the frequency side.  This simplifies our analysis.    Moreover,  we choose this  perturbation  to be real valued and symmetric on the Fourier side i.e.  $\F\vec{\phi}_{0,N}(-\cdot)=\F\vec{\phi}_{0,N}$.  This enables us to consider more general nonlinearity $u^\rho(\bar{u})^{\sigma-\rho}$ as  compared to $u^\sigma$ in \cite{forlano2019remark}.
   In order to get infinite loss of regularity (specifically,  while showing the $X_\theta^{p,q}$-norm of the solution is arbitrarily large for \textit{all} $\theta$),
we restrict the Fourier transform of the solution 
at a particular frequency (say at $n=e_1$).  This will allow us to
compare two discrete Lebesgue norms with \textit{different} weights 
i.e. $$\|\<\cdot\>^\theta f\|_{\ell^q({n=e_1})}=2^{(s-\theta)/2}\n{\<\cdot\>^s f}_{\ell^q(n=e_1)} \quad  \text{for   all}  \ \theta \in \  \mathbb R$$ 
which eventually leads to infinite loss of regularity.
In addition to this, it also allows us to use Plancherel Theorem to achieve the lower bound estimate for the second nontrivial Picard iterate in $W_s^{2,q}$-norm, see Lemma \ref{dt2'}. 
We note that a similar idea was used  by  Kishimoto in \cite[Appendix A]{kishimoto2019remark} in the context of  NLS to achieve norm inflation at zero with infinite loss of regularity in $H^s(\m)$.  Kishimoto  \cite{kishimoto2019remark} used modulation space $M^{2,1}$ to justify the convergence  power series expansion while in the present paper we use Wiener algebra $\mathcal{F}L^1$ as in \cite{OhT, forlano2019remark}.
We employ these ideas   together with the refinement  of ideas used in  \cite{forlano2019remark, OhT, kishimoto2019remark, bhimani2020norm} and properties of $\widehat{w}_s^{p,q}, W^{p,q}_s$ to prove Theorem \ref{dst}.

 Given $\lambda>0,$ if $u$ solves  \eqref{nlw}, then scaling  $u_{\lambda}(t,x):=\lambda^{\frac{2}{\sigma-1}}u(\lambda t, \lambda x)$ also solves \eqref{nlw} with rescaled  initial data $\lambda^{\frac{2}{\sigma-1}} (u_0(\lambda x), u_1(\lambda x)).$  This scaling leaves the homogeneous Sobolev $\dot{H}^s$ space invariant when  $s=s_c= d/2- 2/(\sigma-1).$ The  ill-posedness below the scaling critical regularity $s_c$ has been studied in \cite{LHans, Lebeau, BurqR}.  
 We note that NI with finite loss of regularity at zero initial data for \eqref{nlw} is  initiated by G.  Lebeau in \cite{Lebeau}\footnote{This result is stronger than finite loss of regularity at zero as it is obatied with a single datum instead of a sequence of intial data converging to zero.},  see Theorem B (3).   While we initiate,  to the best of the authors' knowledge,  NI  with infinite loss of regularity at general initial data for \eqref{nlw} in the present paper.  Our study of infinite loss of regularity for \eqref{nlw} is inspired form known NLS case in \cite{CarlesJems,  CarlesT,kishimoto2019remark} while the inspiration of NI at general initial data comes from \cite{OhT, forlano2019remark}.  We also note that  the details of proofs Lebeau in \cite{Lebeau} does not seem to work for negative regularity $s<0.$  In \cite{CarlesJems,  bhimani2020norm, CarlesT}, Carles and his collaborators  have  used geometric optics approach to  establish  infinite loss of regularity  for NLS.  Recently in  \cite{bhimani2021strong}, Bhimani and Haque  have used the Fourier analytic approach to  establish infinite loss of regularity for fractional  Hartree and cubic NLS for some negative  regularity (i.e for all $s<-\epsilon<0$ for some $\epsilon>0$).   On the other hand Theorem \ref{dst} establishes infinite loss of regularity for NLW \eqref{nlw} for \textit{all} $s<0.$  This is in  strict contrast as compared to known NLS case results  proved in \cite{bhimani2021strong,bhimani2020norm,CarlesJems}   and  Theorem \ref{dst} thus reveal new phenomena for NLW \eqref{nlw}.  We conclude   our discussion  with the following remarks.
\begin{remark} In order to get the upper bound  for each Picard iterate  we must have $s \leq 0$ (Lemma \ref{d1}) while  we get the lower bound for one dominating Picard iterate for all $s\in \mathbb R$ (Lemma \ref{d2'}).  The restriction $s<0$  will be required at the final stage in the proof of Theorem \ref{dst}.  Our  approach  thus may not work for  $s\geq 0$. 
\end{remark}
\begin{remark}In \cite[Theorem 1.6]{bhimani2020norm}, Bhimani and Carles have established infinite loss of regularity, via geometric optics approach,  for NLS in $M^{p,q}_s$ for all $1\leq p, q \leq \infty$ and for some $s<0.$   This  somewhat  indicates that restriction on $p=2$ in $M_s^{p,q}, W^{p,q}_s$  in Theorem \ref{dst} is\textit{ just} due to our approach and we believe that  Theorem \ref{dst} is also true for any $1\leq p \leq \infty$ in $M_s^{p,q}, W^{p,q}_s$.  In fact, by taking $p=2,$ we can work on the Fourier side due to  Plancherel theorem.  This makes our analysis somewhat  simple.  On the other hand,  for $p\neq 2,$  we do not know how  to deal with frequency-uniform decomposition operators $\oldsquare_n$.  Also taking Theorem A into account, the case   $s=0$ with $q'\leq(\sigma-1) d$  in Theorem \ref{dst}  remains open.   We plan to address  these issues in our future work.
\end{remark}

\begin{remark} Our method proof should be applicable to \eqref{nlw} on special domain $Z=\mathbb R^{d_1} \times \mathbb T^{d_2} \ (d_1, d_2 \in \mathbb N \cup \{0\})$  with more general type non-linearity  $\sum_{j=1}^n \nu_ju^{\rho_j}(\bar{u})^{\sigma_j-\rho_j}$ where $\nu_j\in\C$, $\sigma_j \in \mathbb N$, $\rho_j\in\N\cup\{0\}$ with $\sigma_j\geq\max(\rho_j,2)$; as in the NLS case\cite{kishimoto2019remark}.
\end{remark}
\begin{remark}\label{stft}
The  STFT  of a  $f\in \mathcal{S}'(\m)$ with
 respect to a window function $0\neq g \in {\mathcal S}(\m)$ is defined by 
\begin{equation*}\label{stft}
V_{g}f(x,y)= \int_{\m} f(t) \overline{T_xg(t)} e^{- 2\pi i y\cdot t}dt,  \  (x, y) \in \m \times \widehat{\m}
\end{equation*}
 whenever the integral exists.  Here, $T_xg(t)=g(tx^{-1})$ is the translation operator on $\m$.  It is known \cite[Proposition 2.1]{wang2007global}, \cite{feichtinger1983modulation} that
\[ \|f\|_{M^{p,q}_s}\asymp  \left\| \|V_gf(x,y)\|_{L^p(\m)} \langle y \rangle^s \right\|_{L^q(\widehat{\m})} \quad \text{and} \quad  \|f\|_{W^{p,q}_s(\m)}\asymp \left\| \|V_gf(x,y) \langle y \rangle^s\|_{L^q(\widehat{\m})}  \right\|_{L^p(\m)}. \] 
 The definition of the modulation space  is independent of the choice of 
the particular window function, see  \cite[Proposition 11.3.2(c)]{grochenig2013foundations}.
\end{remark}
\section{Key Lemmas}
For $u_0,u_1\in \mathcal{S}(\mathbb R^{d}),$ (or $\in C^\infty(\td)$ in torus case) the wave propagator
$S(t)(u_0,u_1)$ is given by $S(t)(u_0,u_1)=\cos(t|\nabla|)u_0+\frac{\sin(t|\nabla|)}{|\nabla|}u_1$, in other words
$$\F S(t)(u_0,u_1)(\xi)=\cos(t|\xi|)\F u_0(\xi)+\frac{\sin(t|\xi|)}{|\xi|}\F u_1(\xi),\quad\xi\in\mh,\  t\in \mathbb R.$$
Let  $\mathcal{X}_s^{p,q}(\m)=\widehat{w}_s^{p,q}(\m)\times{\widehat{w}_{s-1}^{p,q}}(\m)$ with the norm
\[\n{\overrightarrow{u_0}}_{\mathcal{X}_s^{p,q}}=\n{u_0}_{{\widehat{w}_s^{p,q}}}+\sqrt{2}\n{u_1}_{{\widehat{w}_{s-1}^{p,q}}} \quad (1\leq p,q \leq \infty,  s\in \mathbb R).
\]
When $s=0$, we  write $\mathcal{X}^{p,q}(\m)=\mathcal{X}_0^{p,q}(\m)$.

\begin{lemma}\label{w1}
Let $\mathcal{X}_s^{p,q}(\m)$ be defined as above and $0\leq t\leq1.$  Then 
$\n{S(t)(\vec{u}_0)}_{\widehat{w}_s^{p,q}}\leq\n{\vec{u}_0}_{\mathcal{X}_s^{p,q}}$.
\end{lemma}  
\begin{proof}
Since  $\abs{\cos(t|\xi|)}\leq1$ and 
for $0\leq t\leq1,$
\[
\frac{\abs{\sin(t|\xi|)}}{|\xi|}(1+|\xi|^2)^{1/2}\leq\begin{cases}
t(1+|\xi|^2)^{1/2}\leq \sqrt{2}& \text{if }|\xi|\leq1\\
\frac{(1+|\xi|^2)^{1/2}}{|\xi|}\leq \sqrt{2}& \text{if }|\xi|\geq1
\end{cases}
\] the result follows from the definition of $\widehat{w}^{p,q}_s(\m)$.
\end{proof}

For $f_1,\cdots,f_{\sigma}\in\mathcal{S}(\rd)$ (or $\in C^\infty(\td)$ in torus case), we define the multilinear operator $\h_\sigma$ associated to the nonlinearity in \eqref{nlw} as follows
\[ \h_{\sigma} (f_1,\dots,f_{2\sigma +1})= \prod_{\ell=1}^{\rho} f_{\ell} \prod_{m=\rho+1}^{\sigma} \bar{f}_m.\]
When $f_1=\cdots= f_{\sigma}=f$, we write $\h_\sigma(f_1,\dots,f_\sigma)=\h_\sigma(f)$. 
We set 
\[
\mathcal{N}_\sigma(u_1,\dots,u_\sigma)(t)=\int_0^t\frac{\sin((t-\tau)|\nabla|)}{|\nabla|}\h_{\sigma}(u_1(\tau),\cdots,u_\sigma(\tau)) d\tau
\]
and write $\mathcal{N}_\sigma(u_1,\dots,u_\sigma )=\mathcal{N}_\sigma(u)$ for $u=u_1=\cdots=u_\sigma$.
Recall that solution of \eqref{nlw} satisfies
\begin{equation}\label{int}
u(t)=S(t)(u_0,u_1)\pm \int_0^t\frac{\sin((t-\tau)|\nabla|)}{|\nabla|}\h_{\sigma}(u(\tau))d\tau=S(t)(u_0,u_1)\pm\mathcal{N}_\sigma(u)(t)
\end{equation}


\begin{definition}[Picard iteration]\label{imd} 
Let us set
$S_1[\vec{u}_0](t)= S(t)(\vec{u}_0)$ 
and 
\[ S_k[\vec{u}_0](t) = \sum_{{k_1, \dots, k_\sigma \geq 1}\atop{k_1+\dots+k_\sigma =k}}\int_0^t \frac{\sin((t-\tau)|\nabla|)}{|\nabla|}\h_{\sigma}\left(  S_{k_1}[\vec{u}_0],...,S_{k_\sigma}[\vec{u}_0] \right)(\tau) d\tau \quad (k\geq 2).\]
\end{definition}

\begin{remark}
The empty sums in Definition \ref{imd} are considered as zeros. In view of this one can see that $S_{(\sigma-1)\ell+2}[\vec{u}_0]=S_{(\sigma-1)\ell+3}[\vec{u}_0]=\cdots=S_{(\sigma-1)(\ell+1)}[\vec{u}_0]=0$ for all $\ell\in\mathbb{N}\cup\set{0}$. 
\end{remark}

\begin{lemma}[Algebra property]\label{algebra}
The spaces $\widehat{w}^{p,q}_s(\m)$ is a pointwise $\F L^1$-module with norm inequality
$$ \n{fg}_{\widehat{w}^{p,q}_s}\leq \n{f}_{\F L^1}\n{g}_{\widehat{w}^{p,q}_s} \quad  (1\leq p, q \leq \infty, s\in \mathbb R).$$
In particular, $\F L^1$ is an algebra under pointwise multiplication, i.e. $\n{fg}_{\F L^1}\leq \n{f}_{\F L^1}\n{g}_{\F L^1}$. 
\end{lemma}
\begin{proof}
Follows from Young's inequality.
\end{proof}

\begin{lemma}[See \cite{kishimoto2019remark}] \label{DS}
  \label{it0}  Let $\{b_k\}_{k=1}^{\infty}$ be a sequence of nonnegative  real numbers such that 
\[ b_k \leq C  \sum_{{k_1,\dots, k_{\sigma} \geq
      1}\atop{k_1+\dots+ k_{\sigma} =k}}  b_{k_1}\cdots
  b_{k_{\sigma}} \quad \forall\ k \geq 2.\]
Then  we have 
$b_k \leq b_1 C_0^{k-1}$, for all $k \geq 1$, where $C_0=
  \frac{\pi^2}{6} \(C \sigma^2\)^{1/(\sigma-1)} b_1$. 
\end{lemma}
\begin{lemma}\label{w2}
For $0\leq t\leq1$, for $k\geq\sigma$ one has
$\n{S_k[\vec{u}_0](t)}_{\widehat{w}_s^{p,q}} \leq C^{k}t^{2\frac{k-1}{\sigma-1}}\n{\vec{u}_0}_{\overrightarrow{\F L}^1}^{k-1}\n{\vec{u}_0}_{\mathcal{X}_s^{p,q}}.$
\end{lemma}  
\begin{proof}
Let $\{a_k\}$ be a  sequence of nonnegative  real numbers such that 
\[ a_1=1,\quad a_k = \frac{\sigma-1 }{2k-\sigma-1}  \sum_{{k_1,\dots, k_\sigma \geq
      1}\atop{k_1+\dots+ k_\sigma=k}}  a_{k_1}\cdots
  a_{k_\sigma} \quad \forall\ k \geq 2\] where $C>1$ to be chosen later.
  By  Lemma \ref{DS} \eqref{it0}, we have $a_k\leq c_0^{k}$ for some $c_0=c_0(\sigma)>0$. 
In view of this it is enough to prove:
$\n{S_k[\vec{u}_0](t)}_{\widehat{w}_s^{p,q}} \leq a_kt^{2\frac{k-1}{\sigma-1}}\n{\vec{u}_0}_{\overrightarrow{\F L}^1}^{k-1}\n{\vec{u}_0}_{\mathcal{X}_s^{p,q}}. $
By Definition \ref{imd} and the fact $\abs{\sin \tau}\leq|\tau|$ together with Lemma \ref{algebra}, 
we have 
\begin{align*}
&\n{ S_k[\vec{u}_0](t)}_{\widehat{w}_s^{p,q}}
\leq \sum_{{k_1, \dots, k_\sigma \geq 1}\atop{k_1+\dots+k_\sigma =k}}\int_0^t \abs{t-\tau}\n{ S_{k_1}[\vec{u}_0](\tau)}_{\widehat{w}_s^{p,q}}\prod_{\ell=2}^\sigma\n{S_{k_\ell}[\vec{u}_0](\tau)}_{\F L^1}d\tau.
\end{align*}
Therefore, by Lemma \ref{w1}, we have 
\begin{align*}
\n{ S_\sigma[\vec{u}_0](t)}_{\widehat{w}_s^{p,q}}\leq t\int_0^t\n{S_1[\vec{u}_0](\tau)}_{\widehat{w}_s^{p,q}}\prod_{\ell=2}^\sigma\n{S_1[\vec{u}_0](\tau)}_{\F L^1}d\tau\leq t^2\n{\vec{u}_0}_{\overrightarrow{\F L}^1}^{\sigma-1}\n{\vec{u}_0}_{\mathcal{X}_s^{p,q}}.
\end{align*}
Since $a_\sigma=1$,  the claim is true for $k=\sigma$. Assume that the claim  is true upto the level $k-1$. Then
\begin{align*}
\n{ S_k[\vec{u}_0](t)}_{\widehat{w}_s^{p,q}}
&\leq \sum_{{k_1,\dots, k_\sigma  \geq
      1}\atop{k_1+\dots+ k_\sigma  =k}} \n{\vec{u}_0}_{\overrightarrow{\F L}^1}^{k_1-1}\n{\vec{u}_0}_{\mathcal{X}_s^{p,q}}\prod_{\ell=2}^\sigma\n{\vec{u}_0}_{\overrightarrow{\F L}^1}^{k_\ell}
      a_{k_1}\cdots a_{k_\sigma}t\int_0^t\tau^{2\frac{k-\sigma}{\sigma-1}}d\tau \\
      &= a_k\n{\vec{u}_0}_{\overrightarrow{\F L}^1}^{k-1}\n{\vec{u}_0}_{\mathcal{X}_s^{p,q}}t^{2\frac{k-1}{\sigma-1}}.
\end{align*}
Hence,  the claim is true at the level $k$. This completes the proof.
\end{proof}

\begin{lemma}\label{uc}
 If $0<T\ll \min(1,M^{-\frac{\sigma-1}{2}}),$ then  for any  $\vec{u_0} \in \overrightarrow{\F L}^1$ with  $\|\vec{u}_0\|_{\overrightarrow{\F L}^1} \leq M$, there exists a unique solution $u$ to  integral equation \eqref{int} given by
\begin{equation}
  \label{picard-sum}
  u= \sum_{k=1}^{\infty} S_k[\vec{u}_0]= \sum_{\ell =0}^{\infty} S_{2\sigma \ell +1}[\vec{u}_0] 
\end{equation}
which converges absolutely in $C([0, T], \F L^1).$
\end{lemma}
\begin{proof}
The  proof goes in a similar line as the proof of Lemma 2.4 in \cite{  forlano2019remark}.  Since the nonlinearity is different in our case,  we shall briefly present the proof for the convenience of reader. \\
Define 
\[\Psi(u)(t)=S(t)[\vec{u}_0]\pm\int_0^t\frac{\sin((t-\tau)|\nabla|)}{|\nabla|}\h_\sigma(u(\tau))d\tau.\]
Let $0< T\leq1.$ By Lemma \ref{w1} and following the proof of Lemma \ref{w2}, we obtain
\begin{align*}
\n{\Psi(u)}_{C([0,T],\F L^1)}&\lesssim \n{\vec{u}_0}_{\overrightarrow{\F L}^1}+T^2\n{u}_{C([0,T],\F L^1)}^{\sigma},\\
\n{\Psi(u)-\Psi(v)}_{C([0,T],\F L^1)}&\lesssim T^2\max\(\n{u}_{C([0,T],\F L^1)}^{\sigma-1},\n{v}_{C([0,T],\F L^1)}^{\sigma-1}\)\n{u-v}_{C([0,T],\F L^1)}.
\end{align*} Then considering the ball
$B_{2M}^T=\set{\phi\in C([0,T],\F L^1):\n{\phi}_{C([0,T],\F L^1)}\leq 2M}$ with  $T^2,T^2M^{\sigma-1}\ll1$, we find a unique fixed point of $\Psi$ in $B_{2M}^T.$  Hence, the solution to \eqref{int}. This proves the existence of unique solution.  It is not hard to show that this solution is given by \eqref{picard-sum} (see for e.g. \cite[Lemma 2.4]{forlano2019remark} ). This completes the proof.
\end{proof}
\section{The proof of  Theorem \ref{dst}}
We first prove NI with infinite loss of regularity at general data in $\overrightarrow{\F L^1}(\m)\cap\mathcal{X}_s^{p,q}(\m)$. Subsequently, for general data in $\mathcal{X}_s^{p,q}(\m)$ we use the density of $ \overrightarrow{\F L^1}(\m)\cap\mathcal{X}_s^{p,q}(\m)$ in $\mathcal{X}_s^{p,q}(\m)$ ($s<0$). So let us begin with $\vec{u}_0\in\overrightarrow{\F L^1}(\m)\cap\mathcal{X}_s^{p,q}(\m)$. 
Let $N,R\gg1$, $Q=[-1,1]^d$ and $e_1=(1,0,...,0)\in \rd. $ Set
$\Sigma_N=\{\pm Ne_1,\pm2Ne_1\}$
and  \begin{equation}\label{kid}
\F\phi_{0,N}=R\chi_{\Omega_N}\quad\text{with}\quad\Omega_N=\bigcup_{\eta\in\Sigma_N}(\eta+Q).
\end{equation}
Note that $\F{\phi_{0,N}}=\F{\phi_{0,N}}(-\cdot)$ by the symmetry of the set $Q$ and $\Sigma$ and
\begin{align}\label{1}
\|\phi_{0,N}\|_{\widehat{w}_s^{p,q}} =R  \|\|\chi_{n+Q_1} \widehat{\phi}_{0,N}\|_{L^p} \<n\>^s\|_{\ell^q}
= R \(\sum_{n\in\mathbb Z^d} \left| (n+Q_1)\cap \Omega_N \right|^{\frac{q}{p}} \<n\>^{qs}\)^{1/q} \sim  R N^s.
\end{align}
We take \begin{equation}\label{data}
\vec{\phi}_{0,N}=(\phi_{0,N},0)\quad\text{and}\quad \vec{u}_{0,N}=\vec{u}_0+\vec{\phi}_{0,N}.  
\end{equation}
\begin{lemma}[See Lemma 3.6 in \cite{kishimoto2019remark}]\label{kms}
Let $\vec{\phi}_{0,N}$ be given by \eqref{data}. Then there exists $C>0$ such that for 
all $k\in\N$, we have 
\[\left|  \operatorname{supp}\F{ S_k[\vec{\phi}_{0,N}]} (t)\right| \leq C^k
   ,\quad \forall t\geq 0.\]
\end{lemma}
\begin{proof}
Note that $\operatorname{supp}\F{ S_1[\vec{\phi}_{0,N}]} (t)\subset\operatorname{supp}\F u_0;$ which is contained in $4$ cubes with volumes $2^d$.  Hence, $\left|  \operatorname{supp}\F{ S_1[\vec{\phi}_{0,N}]} (t)\right| \leq 2^d4 ,\quad \forall t\geq 0$. Thus, it is enough to prove the following claim: $\operatorname{supp}\F{ S_k[\vec{\phi}_{0,N}]} (t)$ is contained in $4c_d^{k-1}$ number of cubes with volume $2^d$. 
Clearly the claim is true for $k=1$. Assume that the claim is true upto $k-1$ stage. Then
\begin{align*}\label{6}
\operatorname{supp}\F{ S_k[\vec{\phi}_{0,N}]} (t)\subset\sum_{{k_1,\dots, k_{2\sigma +1} \geq
      1}\atop{k_1+\dots+ k_{2\sigma +1} =k}} \operatorname{supp}v_{k_j}(t)
\end{align*}where $v_{k_{\ell}}$ is either $\F{S_{k_{\ell}}[\vec{\phi}_{0,N}]}$
or $\F{\overline{ S_{k_{\ell}}[\vec{\phi}_{0,N}]}}$. Using induction we conclude that the set in RHS is contained in \begin{align*}
d^{2\sigma}\prod_{{k_1,\dots, k_{2\sigma +1} \geq
      1}\atop{k_1+\dots+ k_{2\sigma +1} =k}} 4c_d^{k_j-1}=4^{2\sigma+1}d^{2\sigma}c_d^{k-2\sigma-1}=4(4d)^{2\sigma}c_d^{k-2\sigma-1}
\end{align*}number of cubes with volume $2^d$. Set $c_d=4d$ and $C=c_d=4d$ to conclude.
\end{proof}

\subsection{Estimates in $\widehat{w}_s^{p,q}$}
The next result is the analogue of \cite[Lemma~3.7]{kishimoto2019remark}.
\begin{lemma} \label{d1}Let $\vec{u}_{0,N}$ be given by \eqref{data}, 
 $s\leq0$,
  $1\leq p,q \leq  \infty$ and $0 \leq t \leq 1$. Then there exists $C>0$ independent of $R,N,t$ 
   such that followings hold:
      \begin{enumerate}
      \item $\n{\vec{u}_{0,N}-\vec{u}_0}_{\mathcal{X}_s^{p,q}}\lesssim RN^s$\label{1d1}
      \item $\|S_1[\vec{u}_{0,N}](t)\|_{\widehat{w}_s^{p,q}} \lesssim  1+RN^s$\label{2d1}
      \item $\|S_{\s}[\vec{u}_{0,N}](t)-S_{\s}[\vec{\phi}_{0,N}](t)\|_{\widehat{w}_s^{p,q}}\lesssim t^2R^{\sigma-1}$\label{3d1}
      \item $\|S_k[\vec{u}_{0,N}](t)\|_{\widehat{w}_s^{p,q}} \lesssim C^kR^{k}t^{2\frac{k-1}{\sigma-1}}$.\label{4d1}
      \end{enumerate}
\end{lemma}
\begin{proof}
\eqref{1d1} follows from \eqref{1}.
By Lemma \ref{w1} and \eqref{1}  we have $\|S_1[\vec{\phi}_{0,N}](t)\|_{\widehat{w}_s^{p,q}}
\lesssim  RN^s.$ Then \eqref{2d1} follows by triangular inequality.
By Lemma \ref{w2} (with $p=q=\infty$, $s=0$) and \eqref{1}, we obtain 
\begin{align*}
\|{S_k[\vec{\phi}_{0,N}](t)}\|_{\widehat{w}_s^{p,q}}
&\leq \sup_{\xi \in \mh} |  \F{S_k[\vec{\phi}_{0,N}]} (t, \xi)| |\operatorname{supp}\F{ S_k[\vec{\phi}_{0,N}]} (t)|_{\mu_{\mh}}^{1/p}\n{\<n\>^s}_{\ell^q((n+Q_1)\cap{\rm supp }\ \F S_k[\vec{u}_0](t))}\\
& \lesssim C^{k}    R^{k} t^{2\frac{k-1}{\sigma-1}}\n{\<n\>^s}_{\ell^q((n+Q_1)\cap{\rm supp }\ \F S_k[\vec{u}_0](t))}.
\end{align*}
where $|A|_{\mu_{\mh}}$ denotes the $\mh$-measure of the set $A$.  Since $s\leq0$, for any bounded set $D\subset \zd$, 
we have 
$\| \langle n\rangle^s \|_{\ell^q (  n \in D )}\leq \| \langle n\rangle^s \|_{\ell^q ( n \in B_D )}$
where $B_D\subset \R^d$ is the minimal ball centred at the origin with $|D|\leq|B_D|.$  
In view of this and Lemma \ref{kms}, 
we obtain 
\begin{equation*}
\| \langle n\rangle^s \|_{\ell^q( \operatorname{supp}
  \widehat{S_k}[\vec{u}_0](t))}\leq
\| \langle n\rangle^s \|_{\ell^q ( \{ |n| \leq C^{k/d} \})}
\lesssim   C^{k/q}. 
\end{equation*}
Therefore 
\begin{align}\label{5d1}
\|{S_k[\vec{\phi}_{0,N}](t)}\|_{\widehat{w}_s^{p,q}}\leq C^kR^{k} t^{2\frac{k-1}{\sigma-1}}.
\end{align}
Now observe that 
\begin{align*}
I_{k}(t):=&S_k[\vec{u}_{0,N}](t)-S_k[\vec{\phi}_{0,N}](t)\\
&=\sum_{{k_1, \dots, k_{\sigma} \geq 1}\atop{k_1+\dots+k_{\sigma} =k}}\mathcal{N}(S_{k_1}[\vec{u}_0+\vec{\phi}_{0,N}],\cdots,S_{k_{\sigma}}[\vec{u}_0+\vec{\phi}_{0,N}])-\mathcal{N}(S_{k_1}[\vec{\phi}_{0,N}],\cdots,S_{k_{\sigma}}[\vec{\phi}_{0,N}])\\
&=\sum_{{k_1, \dots, k_{\sigma} \geq 1}\atop{k_1+\dots+k_{\sigma} =k}}\sum_{(\vec{\psi_1},\cdots,\vec{\psi}_{\sigma})\in\mathcal{C}}\mathcal{N}(S_{k_1}[\vec{\psi}_1],\cdots,S_{k_{\sigma}}[\vec{\psi}_{\sigma}])
\end{align*}
where $\mathcal{C}=\{\vec{u}_0,\vec{\phi}_{0,N}\}^{\sigma}\setminus\{(\vec{\phi}_{0,N},\cdots,\vec{\phi}_{0,N})\}$. Observe that $\mathcal{C}$ has atleast one coordinate as $\vec{u}_0$.
Using  Lemma \ref{w2}
 it follows that
\begin{align*}
\|I_k(t)\|_{\widehat{w}_s^{p,q}}
& \lesssim \sum_{{k_1, \dots, k_{\sigma} \geq 1}\atop{k_1+\dots+k_{\sigma} =k}}t\int_0^t\sum_{(\vec{v}_1,\cdots,\vec{v}_{\sigma})\in \mathcal{C}}\|S_{k_1}[\vec{v}_1]\|_{\widehat{w}_s^{p,q}}\prod_{j=2}^{\sigma}\|S_{k_j}[\vec{v}_j]\|_{{\F L^1}}\\
&\leq (2^{\sigma}-1)2\|\vec{u}_0\|_{\mathcal{X}_s^{p,q}} (\|\vec{u}_0\|_{\overrightarrow{\F L^1}}^{k-1}+\|\vec{\phi}_{0,N}\|_{\overrightarrow{\F L^1}}^{k-1})t\int_0^t\tau^{2\frac{k-\sigma}{\sigma-1}}d\tau\sum_{{k_1, \dots, k_{\sigma} \geq 1}\atop{k_1+\dots+k_{\sigma} =k}}a_{k_1}\cdots a_{k_{\sigma}}\\
&\leq 2^{\sigma+2}a_kt^{2\frac{k-1}{\sigma-1}}R^{k-1}\|\vec{u}_0\|_{\mathcal{X}_s^{p,q}}\leq C^kt^{2\frac{k-1}{\sigma-1}}R^{k-1}\|\vec{u}_0\|_{\mathcal{X}_s^{p,q}}
\end{align*}
as $R\gg1$. Note that \eqref{3d1} is the particular case $k=\sigma$ and \eqref{4d1} follows using the above and \eqref{5d1}. 
\end{proof}

In the proof of  next the lemma we follow the strategy  of  Proposition 3.4 in \cite{forlano2019remark}. Although Proposition 3.4 in \cite{forlano2019remark} considers a different nonlinearity ($u^{\sigma}$), the  symmetry about the origin of the real valued function $\F u_0$ allows us to cover our choice of nonlinearity ($u^{\rho}\bar{u}^{\sigma-\rho}$). We have presented the proof in detail as it will  be used in the proof of the similar estimates in $W_s^{2,q}(\rd)$ spaces (Lemma \ref{dt2'}).
\begin{lemma}\label{d2'}Let $\vec{\phi}_{0,N}$ be given by \eqref{data}, $s\in \mathbb R$,  $1\leq p,q \leq  \infty$, $N^{-1/2}<T\ll 1.$ Then we   have 
\[ \|S_\sigma [\vec{\phi}_{0,N}] (T)\|_{\widehat{w}_s^{p,q}}\geq \left\| \left\lVert \chi_{n+Q_1}(\xi)\F S_\sigma [\vec{\phi}_{0,N}] (T)(\xi)\right\rVert_{L_\xi^p} \<n\>^{s} \right\|_{\ell^q(n=e_1)} 
\gtrsim  R^{\sigma}T^2. \]
\end{lemma}
\begin{proof} We shall first briefly gives the guideline of the proof.
In order to establish the required lower estimate for $\|S_{\sigma} [\vec{\phi}_{0,N}] (T)\|_{\widehat{w}_s^{p,q}}$, we first write $\F S_{\sigma} [\vec{\phi}_{0,N}] (T)(\xi)$ (with $\xi\in Q_2$) in terms of a double sum (see  \eqref{id1}).  We shall naturally arrive summation over $\mathcal{A}$ by the choice of data \eqref{data} with large $N.$ And then over $\mathcal{B}$ by applying  suitable trigonometric identities for cosine functions.  Further we will divide the terms under summation into two category: one collection of good terms (say $I_0$) which helps us to get the lower bound and the other collection  of bad terms ($I_1$) which none the less will have some upper bound. Subsequently  we will choose the time $T$ so that the good terms dominate over the bad terms to achieve the required estimate.

We now  produce the details of the  proof. Note that $\overline{\F u_0}(-\cdot)=\F u_0$. 
Set $$\Gamma_\xi=\left\{(\xi_1,\cdots,\xi_{2\sigma+1})\in\R^{\sigma d}:\sum_{\ell=1}^{\sigma}\xi_\ell=\xi\right\}$$ and $d\Gamma_{\xi}$ denote the $(\sigma-1)$-dimensional Lebesgue measure on the hyperplane $\Gamma_\xi$.
Note that for $\xi\in Q_2$, using $N\gg1$ we have
\begin{align}\label{id1}
&\F S_{\sigma} [\vec{\phi}_{0,N}] (T)(\xi) \nonumber \\
&=\int_0^T\frac{\sin((T-t)|\xi|)}{|\xi|}\left[\(
\Conv^{\rho}_{\ell=1}\cos(t|\cdot|)\F u_0\)\Conv\(
\Conv^{\sigma}_{m=\rho+1}\cos(t|\cdot|)\overline{\F u_0}(-\cdot)\)\right](\xi)dt \nonumber \\
&=\int_0^T\frac{\sin((T-t)|\xi|)}{|\xi|}\int_{\Gamma_\xi}
\prod_{\ell=1}^{\sigma}\cos(t|\xi_\ell|)\F u_0(\xi_\ell)d\Gamma_\xi dt \nonumber \\
&=R^{\sigma}\sum_{\mathcal{A}}\int_0^T\frac{\sin((T-t)|\xi|)}{|\xi|}\int_{\Gamma_\xi}
\prod_{\ell=1}^{\sigma}\cos(t|\xi_\ell|)\1_{\eta_\ell+Q_2}(\xi_\ell)d\Gamma_\xi dt \nonumber \\
&=\frac{R^{\sigma}}{4^{\sigma}}\sum_{\mathcal{A}}\sum_{\mathcal{B}}\int_0^T\frac{\sin((T-t)|\xi|)}{|\xi|}\int_{\Gamma_\xi}
\cos\(t\sum_{\ell=1}^{\sigma}\epsilon_\ell|\xi_\ell|\)
\prod_{\ell=1}^{\sigma}\1_{\eta_\ell+Q_2}(\xi_\ell)d\Gamma_\xi dt
\end{align}
where the sums are taken over 
\begin{align*}
\mathcal{A}=\left\{(\eta_1,\cdots,\eta_{\sigma})\in\Sigma^{\sigma}:\sum_{\ell=1}^{\sigma}\eta_\ell=0\right\},\quad\mathcal{B}=\left\{(\eps_1,\cdots,\eps_{\sigma})\in\{\pm1\}^{\sigma}\right\}
\end{align*}respectively. 
For $\eta=(\eta_1,\cdots,\eta_{\sigma})\in\Sigma^{\sigma}$ set  $$\mathcal{B}_0(\eta)=\left\{(\eps_1,\cdots,\eps_{\sigma})\in\{\pm1\}^{\sigma}:\sum_{\ell=1}^{\sigma}\eps_\ell|\eta_\ell|=0\right\},\quad\mathcal{B}_1(\eta)=\mathcal{B}\setminus\mathcal{B}_0(\eta).$$
Then splitting the inner sum in  \eqref{id1} over $\mathcal{B}_0(\eta)$ and $\mathcal{B}_1(\eta)$ we write
\begin{align*}
\F S_{\sigma} [\vec{\phi}_{0,N}] (T)(\xi)=\frac{R^{\sigma}}{4^{\sigma}}\sum_{\mathcal{A}}\(I_0(\eta,T,\xi)+I_1(\eta,T,\xi)\).
\end{align*}
Note that for each $\eta\in\mathcal{A}$ the set $\mathcal{B}_0(\eta)$ is non empty. This is because for $\eta\in\mathcal{A}$, $\sum_{\ell=1}^{\sigma}\eta_\ell=0$ and then $\(\frac{(\eta_1)_1}{|\eta_1|},\cdots,\frac{(\eta_\sigma)_1}{|\eta_\sigma|}\)\in\mathcal{B}_0(\eta)$ (here $(\eta_\ell)_1$ denotes the first coordinate of $\eta_\ell$).  For a fixed $\eps=(\eps_1,\cdots,\eps_{\sigma})\in\mathcal{B}_0(\eta)$ and $\xi_\ell\in\eta_\ell+Q_2$ one has (using triangular inequality)
\[\abs{\sum_{\ell=1}^{\sigma}\eps_\ell|\xi_\ell|}\leq\abs{\sum_{\ell=1}^{\sigma}\eps_\ell|\eta_\ell|}+\abs{\sum_{\ell=1}^{\sigma}|\xi_\ell-\eta_\ell|}\lesssim 1.\]
Therefore, for $0\leq t\leq T\ll 1$, we have $\cos\(t\sum_{\ell=1}^{\sigma}\eps_\ell|\xi_\ell|\)\geq\frac{1}{2}$.
On the other hand, we have 
$$\frac{\sin((T-t)|\xi|)}{|\xi|}\gtrsim T-t$$ provided $0\leq t<T\ll 1$ and $\xi\in Q_2$. Hence,  for $\xi\in Q_2$, 
we have  
\begin{align*}
I_0(\eta,T,\xi)\gtrsim\sum_{\eps\in \mathcal{B}_0}\int_0^T(T-t)\int_{\Gamma_\xi}\prod_{\ell=1}^{\sigma}\1_{\eta_\ell+Q_2}(\xi_\ell)d\Gamma_\xi dt\gtrsim T^2\1_{Q_2}(\xi)
\end{align*} 
as $\1_{\alpha+Q_2}*\1_{\beta+Q_2}\geq c_d\1_{\alpha+\beta+Q_2}$ for $\alpha,\beta\in\rd$. Therefore, for $\xi\in Q_2$, $N\gg1$ and $0<T\ll 1$, we have\begin{equation}\label{8}
\frac{R^{\sigma}}{4^{\sigma}}\sum_{\mathcal{A}}I_0(\eta,T,\xi)\gtrsim T^2R^{\sigma}\1_{Q_2}(\xi).
\end{equation}
Note that for $\eps\in\mathcal{B}_1(\eta)$ and $\xi_\ell\in\eta_\ell+Q_2$ one has $\abs{\sum_{\ell=1}^{\sigma}\eps_\ell|\xi_\ell|}\lesssim N$, this together with 
\[\sum_{\ell=1}^{\sigma}\eps_\ell|\xi_\ell|=\sum_{\ell=1}^{\sigma}\eps_\ell|\eta_\ell|+\sum_{\ell=1}^{\sigma}\eps_\ell(|\xi_\ell|-|\eta_\ell|)=\sum_{\ell=1}^{\sigma}\eps_\ell|\eta_\ell|+\mathcal{O}(1)\]
implies \begin{equation}
\abs{\sum_{\ell=1}^{\sigma}\eps_\ell|\xi_\ell|}\sim N.
\end{equation} Therefore using this  in 
\begin{align*}
& I_1(\eta,T,\xi)\\
&=\sum_{\eps\in \mathcal{B}_1(\eta)}\frac{1}{2|\xi|}\int_0^T\int_{\Gamma_\xi}\left[\sin\left(T|\xi|-t\left(|\xi|+\sum_{\ell=1}^{\sigma}\eps_\ell|\xi_\ell|\right)\right)
+\sin\left(T|\xi|-t\left(|\xi|-\sum_{\ell=1}^{\sigma}\eps_\ell|\xi_\ell|\right)\right)\right]\\&\hspace{1.2cm} 
\times\prod_{\ell=1}^{\sigma}\1_{\eta_\ell+Q_2}(\xi_\ell)d\Gamma_\xi dt\\
&=\sum_{\eps\in \mathcal{B}_1(\eta)}\frac{1}{2|\xi|}\int_{\Gamma_\xi}\left[\frac{\cos\left(T|\xi|-t\left(|\xi|+\sum_{\ell=1}^{\sigma}\eps_\ell|\xi_\ell|\right)\right)}{|\xi|+\sum_{\ell=1}^{\sigma}\eps_\ell|\xi_\ell|}
+\frac{\cos\left(T|\xi|-t\left(|\xi|-\sum_{\ell=1}^{\sigma}\eps_\ell|\xi_\ell|\right)\right)}{|\xi|-\sum_{\ell=1}^{\sigma}\eps_\ell|\xi_\ell|}\right]_{t=0}^T\\&\hspace{1.2cm} 
\times\prod_{\ell=1}^{\sigma}\1_{\eta_\ell+Q_2}(\xi_\ell)d\Gamma_\xi \\
\end{align*} 
we get $|I_1(\eta,T,\xi)|\lesssim N^{-1}\1_{2^{\sigma-1} Q_2}(\xi)$  for $\xi\in \frac{3}{4}e_1+Q_{\frac{1}{2}}$ as $\1_{\alpha+Q_2}*\1_{\beta+Q_2}\leq c_d\1_{\alpha+\beta+2Q_2}$.
Hence for $\xi\in\frac{3}{4}e_1+Q_{\frac{1}{2}},$
\begin{equation}\label{9}
\frac{R^{\sigma}}{4^{\sigma}}\abs{\sum_{\mathcal{A}}I_1(\eta,T,\xi)}\lesssim N^{-1}R^{\sigma}\1_{2^{\sigma-1}  Q_{2}}(\xi).
\end{equation}
Therefore using \eqref{8} we have for $\xi \in\frac{3}{4}e_1+Q_{\frac{1}{2}}\subset Q_2$
\begin{equation}\label{W1}
|\F S_{\sigma} [\vec{\phi}_{0,N}] (T,\xi)|\gtrsim T^2R^{\sigma}
\end{equation}
provided $T^2\gg N^{-1}$ and $0<T\ll 1$.
Thus we conclude 
\begin{align}\label{10}
\|S_{\sigma} [\vec{\phi}_{0,N}] (T)\|_{\widehat{w}_s^{p,q}}\geq\left\| \left\lVert \chi_{n+Q_1}(\xi)\F S_{\sigma} [\vec{\phi}_{0,N}] (T)(\xi)\right\rVert_{L_\xi^p} \<n\>^{s} \right\|_{\ell^q(n=e_1)} 
\gtrsim T^2R^{\sigma}
\end{align}if $N^{-1/2}<T\ll 1$. 
\end{proof}

\subsection{Estimates in $W_s^{2,q}$}
\begin{lemma}[inclusion]\label{INC} Let $p,q,q_1,q_2\in[1,\infty]$ and $s\in\R$. Then (1) $\|f\|_{W_s^{2,q}}\leq \|f\|_{\widehat{w}_s^{2,q}}$ if $q\leq 2$ (2) $\|f\|_{W_s^{p,q_1}}\lesssim\|f\|_{W_s^{p,q_2}}$ if $q_1\geq q_2$.
\end{lemma}
\begin{proof}
(1) is a consequence of Minkowski inequality and Plancherel theorem  whereas (2) follows from the fact that $\ell^{q_2}\hookrightarrow\ell^{q_1}$ if $q_1\geq q_2$. 
\end{proof}
\begin{lemma}\label{dt1'} 
Let $\vec{u}_{0,N}$ be given by \eqref{data}, 
 $s\leq0$,
  $1\leq q \leq  \infty$ and $0 \leq t \leq 1$. Then there exists $C>0$ independent of $R,N,t$ 
   such that followings hold:
      \begin{enumerate}
      \item $\n{\vec{u}_{0,N}-\vec{u}_0}_{{\mathcal W}_s^{2,q}}\lesssim RN^s$, where $\mathcal W_s^{2,q}=W_s^{2,q}\times W_{s-1}^{2,q}$
      \item $\|S_1[\vec{u}_{0,N}](t)\|_{{W}_s^{2,q}} \lesssim  1+RN^s$
      \item $\|S_{\sigma}[\vec{u}_{0,N}](t)-S_{\sigma}[\vec{\phi}_{0,N}](t)\|_{{W}_s^{2,q}}\lesssim t^2R^{\sigma-1}$
      \item $\|S_k[\vec{u}_{0,N}](t)\|_{{W}_s^{2,q}} \lesssim C^kR^{k}t^{2\frac{k-1}{\sigma-1}}$.
      \end{enumerate}
\end{lemma}
\begin{proof}
By Lemma \ref{INC}, we have 
\begin{align*}
\|\vec{u}_{0,N}-\vec{u}_0\|_{W^{2,q}_s} \lesssim\begin{cases}
\|\vec{u}_{0,N}-\vec{u}_0\|_{\widehat{w}^{2,q}_s}\lesssim RN^s&\text{ for }q\in[1,2]\\
 \|\vec{u}_{0,N}-\vec{u}_0\|_{W^{2,2}_s} \lesssim RN^s&\text{ for }q\in(2,\infty]
\end{cases}
\end{align*}  
using Lemma \ref{d1} \eqref{1d1}. Similarly the other estimates also follow from Lemma \ref{d1}.
\end{proof}
\begin{lemma}\label{dt2'}
Let $\vec{\phi}_{0,N}$ be given by \eqref{data}, $s\in \mathbb R$,  $1\leq q \leq  \infty$, $N^{-1/2}<T\ll 1.$ Then we   have 
\begin{align*} 
\|S_{\sigma} [\vec{\phi}_{0,N}] (T)\|_{W^{2,q}_s}\geq\left\| \lVert\oldsquare_n S_{\sigma} [\vec{\phi}_{0,N}] (T)(\xi)(1+|n|)^{s}\rVert_{\ell^q(n=e_1)}  \right\|_{L_\xi^2}  
 \gtrsim R^{\sigma}T^2.
\end{align*}
\end{lemma}
\begin{proof}
Note that using Plancherel theorem and \eqref{W1} we have
\begin{align*}
\|S_{\sigma} [\vec{\phi}_{0,N}] (T)\|_{W^{2,q}_s}&\geq\left\| \lVert\F^{-1} \sigma_n\F S_{\sigma} [\vec{\phi}_{0,N}] (T)(x)(1+|n|)^{s}\rVert_{\ell^q(n=e_1)}  \right\|_{L_x^2} \\
&=2^{s/2}\left\|  \sigma_{e_1}\F S_{\sigma} [\vec{\phi}_{0,N}] (T)(\xi)  \right\|_{L_\xi^2}\gtrsim R^{\sigma }T^2.
\end{align*}
This completes the proof.
\end{proof}

\subsection{Proof of main result}
\begin{proof}[{\bf Proof of Theorem \ref{dst}}]   
We first consider the case when $X_s^{p,q}(\m)=\widehat{w}_s^{p,q}(\m)$. In view of the comment at the beginning of this section, it is enough to prove NI with infinite loss of regularity at $\vec{u}_0 \in \overrightarrow  {\F L^1}(\m)\cap\mathcal{X}_s^{p,q}(\m)$. By
  Lemma~\ref{uc}, we have  the existence of a unique solution  to
  \eqref{nlw} with initial condition given by \eqref{data} in $\F L^1(\m)$ 
   up to time  $T$ whenever  $(\|\vec{u}_0\|_{\overrightarrow  {\F L^1}}+R)^{\frac{\sigma-1}{2}} T \ll 1$ which is implied by (0) $R^{\frac{\sigma-1}{2}} T \ll 1$ as $R\gg1$.  In view
  of  Lemma~\ref{d1} and if $RT^{\frac{2}{\sigma-1}}\ll1$ ($\Leftrightarrow(0)$),  $\sum_{\ell =2}^{\infty}
  \| S_{(\sigma-1) \ell +1}[\vec{u}_{0,N}] (T)\|_{\widehat{w}_s^{p,q}} $ can be
  dominated by the  sum of the geometric series. Specifically, we have  
\begin{align}\label{mv1}
 \sum_{\ell =2}^{\infty}\left\| S_{(\sigma-1) \ell +1}[\vec{u}_{0,N}] (T) \right\|_{\widehat{w}_s^{p,q}}  \lesssim   R  \sum_{\ell=2}^{\infty} (CR)^{(\sigma-1) \ell}T^{2\ell} 
\lesssim  R^{2\sigma-1} T^4.
\end{align}
Note that 
\begin{align*}
\n{u_N(T)}_{\widehat{w}_\theta^{p,q}}\geq\n{\n{\chi_{n+Q_1}\F u_N(T)}_{L^p}\<n\>^\theta}_{\ell^q(n=e_1)}\sim_{\theta,s}\n{\n{\chi_{n+Q_1}\F u_N(T)}_{L^p}\<n\>^s}_{\ell^q(n=e_1)}
\end{align*}
therefore by Lemma \ref{uc} and triangle inequality, we obtain
\begin{align*}
&\n{u_N(T)}_{\widehat{w}_\theta^{p,q}}\\
&\gtrsim \n{\n{\chi_{n+Q_1}\F S_{\sigma}[u_{0,N}](T)}_{L^p}\<n\>^s}_{\ell^q(n=e_1)}-c\bigg(\n{\n{\chi_{n+Q_1}\F S_1[u_{0,N}](T)}_{L^p}\<n\>^s}_{\ell^q(n=e_1)}\\&\ \ \ +\sum_{\ell=2}^\infty\n{\n{\chi_{n+Q_1}\F S_{2\sigma\ell}[u_{0,N}](T)}_{L^p}\<n\>^s}_{\ell^q(n=e_1)}\bigg)\\
&\gtrsim \n{\n{\chi_{n+Q_1}\F S_{\sigma}[u_{0,N}](T)}_{L^p}\<n\>^s}_{\ell^q(n=e_1)}-c\n{S_1 [\vec{u}_{0,N}] (T)}_{\widehat{w}_s^{p,q}}-c\sum_{\ell=2}^\infty\n{S_{(\sigma-1)\ell +1} [\vec{u}_{0,N}] (T)}_{\widehat{w}_s^{p,q}}.
\end{align*}
Assume $m\in\N$ be given.  In order to ensure  
$ \|u(T)\|_{\widehat{w}_\theta^{p,q}} \gtrsim  \| S_{\sigma}[\vec{u}_{0,N}](T)\|_{\widehat{w}_s^{p,q}}  \gg m$
we rely on the conditions 
\begin{numcases}{\n{\n{\chi_{n+Q_1}\F S_{\sigma}[u_{0,N}](T)}_{L^p}\<n\>^s}_{\ell^q(n=e_1)}\gg}
\|S_1[\vec{u}_{0,N}](T)\|_{\widehat{w}_s^{p,q}},\label{rc1}\\
 \Sigma_{\ell =2}^{\infty}\left\| S_{(\sigma-1) \ell +1}[\vec{u}_{0,N}](T) \right\|_{\widehat{w}_s^{p,q}},\label{rc2}\\
m.\label{rc3}
\end{numcases}
To achieve \eqref{rc1}-\eqref{rc3},  we rely on  
 Lemmata \ref{d1}, \ref{d2'}. To use Lemma \ref{d2'}, we impose 
 (i) $N^{-1/2}\ll T\ll 1$. 
In view  of Lemmata \ref{d1}, \ref{d2'}, and \eqref{mv1},   to prove \eqref{rc2}  it is sufficient to have: 
(ii) (a)$R^{\sigma}T^2\gg  R^{2\sigma-1} T^4\Longleftrightarrow   R^{\sigma-1}T^2\ll1 $ ($\Leftrightarrow(0)$) and
(ii) (b) $T^2R^{\sigma-1}\ll  R^{\sigma}T^2\Longleftrightarrow R\gg1$. 
 To achieve \eqref{rc3} we impose
(iii) $R^{\sigma}T^2\gg m $ (along with (ii) (b)).
To ensure $\n{\vec{u}_{0,N}-\vec{u}_0}_{\mathcal{X}_s^{p,q}}<1/m$, in view of \eqref{1} we impose
(iv) $RN^{s}<1/m$.
 At the end (iii)  (along with (ii) (b)) imply \eqref{rc1} using Lemma \ref{d1}.
\\\\ \textbf{{Case $-\frac{1}{\sigma-1}\leq s<0$.}}\\
 Set $R=N^{-s-\delta}, T=N^{\frac{\sigma-1}{2}(s+\frac{\delta}{2})}$ with $0<\delta\ll1$ satisfying $\frac{\sigma+1}{2}\delta<-s$. Note that with $N\gg1$ we have $R\gg1$ and
\begin{align*}
&N^{-\frac{1}{2}}T^{-1}\sim N^{-\frac{1}{2}-\frac{\sigma-1}{2}(s+\frac{\delta}{2})}=N^{-\frac{\sigma-1}{2}(s+\frac{1}{\sigma-1}+\frac{\delta}{2})}\ll1&&\text{ as }s+\frac{1}{\sigma-1}+\frac{\delta}{2}>s+\frac{1}{\sigma-1}\geq0\\
 &T\sim N^{\frac{\sigma-1}{2}(s+{\delta}/{2})}\ll1&&\text{ as }s+\delta/2<-\sigma\delta/2<0\\
& R^{\sigma-1}T^2\sim N^{-(\sigma-1)(s+\delta)+(\sigma-1)(s+\frac{\delta}{2})}=N^{-(\sigma-1)\frac{\delta}{2}}\ll1&&\text{ as } -(\sigma-1)\delta/2<0\\
 &RN^{s}\sim N^{-(s+\delta)+s}=N^{-\delta}\ll\frac{1}{m}&&\text{ as } -\delta<0\\
 &R^{\sigma}T^2\sim N^{-\sigma(s+\delta)+(\sigma-1)(s+\delta/2)}=N^{-s-(\sigma+1)\delta/2}\gg m&&\text{ as }-s-(\sigma+1)\delta/2>0.
\end{align*}
 \textbf{Case}  $s<-\frac{1}{\sigma-1}$.\\
Set $R=N^{\frac{1}{\sigma-1}-\delta}, T=N^{-\frac{1}{2}+\frac{\sigma-1}{4}\delta}$ with $0<\delta\ll1$ satisfying $(\sigma+1)\delta<\frac{2}{\sigma-1}$ and $(\sigma-1)\frac{\delta}{4}<\frac{1}{2}$. Note that with $N\gg1$ we have $R\gg1$ and
\begin{align*}
&N^{-1/2}T^{-1}\sim N^{-1/2+1/2-(\sigma-1)\delta/4}=N^{-(\sigma-1)\delta/4}\ll1&&\text{ as }(\sigma-1)\delta/4>0\\
& T\sim N^{-1/2+(\sigma-1)\delta/4}\ll1&&\text{ as }-1/2+(\sigma-1)\delta/4<0\\
& R^{\sigma-1}T^2\sim N^{1-(\sigma-1)\delta-1+(\sigma-1)\delta/2}=N^{-(\sigma-1)\delta/2}\ll1&&\text{ as } -(\sigma-1)\delta/2<0\\
 &RN^{s}\sim N^{1/(\sigma-1)-\delta+s}\ll\frac{1}{m}&&\text{ as } 1/(\sigma-1)+s<0\\
 &R^{\sigma}T^2\sim N^{\sigma(\frac{1}{\sigma-1}-\delta)-1+(\sigma-1)\frac{\delta}{2}}=N^{\frac{1}{\sigma-1}-(\sigma+1)\frac{\delta}{2}}\gg m&&\text{ as }1/(\sigma-1)-(\sigma+1)\delta/2>0.
\end{align*}
Thus with both the cases the conditions (i)-(iv)  are satisfied and hence we are done with the case $\mathcal{X}_s^{p,q}=\widehat{w}_s^{p,q}$.

For the case $X_s^{p,q}=W_s^{2,q}$ we use same argument as above: 
 Note that using Lemmata \ref{dt1'}, \ref{dt2'}.
\begin{align*}
&\n{u_N(T)}_{W_{\theta}^{2,q}}\\
&\geq\n{\n{ \oldsquare_n u_N(T)\<n\>^\sigma}_{\ell^q(n=e_1)}}_{L^2}\sim_{ \theta,s}\n{\n{ \oldsquare_n u_N(T)\<n\>^s}_{\ell^q(n=e_1)}}_{L^2}\\
&\gtrsim\n{\n{\oldsquare_n S_{\sigma}[u_{0,N}](T)\<n\>^s}_{\ell^q(n=e_1)}}_{L^2}- c\|S_1[u_{0,N}](T)\|_{W_s^{2,q}}- c \sum_{\ell =2}^{\infty}\left\| S_{(\sigma-1)\ell +1}[\vec{u}_{0,N}](T) \right\|_{W_s^{2,q}}\\
&\gtrsim\n{\n{\oldsquare_n S_{\sigma}[u_{0,N}](T)\<n\>^s}_{\ell^q(n=e_1)}}_{L^2}\gg m.
\end{align*}
and $\n{\vec{u}_{0,N}-\vec{u}_0}_{\mathcal{X}_s^{2,q}}<1/m$ provided we choose $R,N,T$ as in the case of $\widehat{w}_s^{p,q}$. 
\end{proof}

\noindent
{\textbf{Acknowledgement}:} D.G. B is thankful to DST-INSPIRE (DST/INSPIRE/04/2016/001507) for the research grant. S. H. acknowledges Dept of Atomic Energy, Govt of India, for the financial support and Harish-Chandra Research Institute for the research facilities provided. Both authors are grateful to professor R\'emi Carles for his thoughtful comments and fruitful discussion on  the topic.  D. G. B. is grateful to professor Tadahiro Oh and Justin Forlano for sending their preprint \cite{JH} and introducing to Fourier amalgam spaces.
\bibliographystyle{siam}
\bibliography{wave}

\begin{thebibliography}{10}

\bibitem{BejenaruTao}
{\sc I.~Bejenaru and T.~Tao}, {\em Sharp well-posedness and ill-posedness
  results for a quadratic non-linear {S}chr\"{o}dinger equation}, J. Funct.
  Anal., 233 (2006), pp.~228--259.

\bibitem{benyi2007unimodular}
{\sc A.~B\'{e}nyi, K.~Gr\"{o}chenig, K.~A. Okoudjou, and L.~G. Rogers}, {\em
  Unimodular {F}ourier multipliers for modulation spaces}, J. Funct. Anal., 246
  (2007), pp.~366--384.

\bibitem{benyi2009local}
{\sc A.~B\'{e}nyi and K.~A. Okoudjou}, {\em Local well-posedness of nonlinear
  dispersive equations on modulation spaces}, Bull. Lond. Math. Soc., 41
  (2009), pp.~549--558.

\bibitem{kassob}
{\sc {\'A}.~B{\'e}nyi and K.~A. Okoudjou}, {\em {M}odulation {S}paces: {W}ith
  {A}pplications to {P}seudodifferential {O}perators and {N}onlinear
  {S}chr\"odinger Equations,}, 2020.

\bibitem{bhimani2020norm}
{\sc D.~G. Bhimani and R.~Carles}, {\em Norm inflation for nonlinear
  {S}chrodinger equations in {F}ourier-{L}ebesgue and modulation spaces of
  negative regularity}, Journal of Fourier Analysis and Applications, 26
  (2020).

\bibitem{bhimani2020hartree}
{\sc D.~G. Bhimani, M.~Grillakis, and K.~A. Okoudjou}, {\em The
  {H}artree--{F}ock equations in modulation spaces}, Communications in Partial
  Differential Equations,  (2020), pp.~1--30.

\bibitem{bhimani2021strong}
{\sc D.~G. Bhimani and S.~Haque}, {\em Strong ill-posedness for fractional
  {H}artree and cubic {NLS} equations}, preprint arXiv:2101.03991,  (2021).

\bibitem{bhimani2016functions}
{\sc D.~G. Bhimani and P.~K. Ratnakumar}, {\em Functions operating on
  modulation spaces and nonlinear dispersive equations}, J. Funct. Anal., 270
  (2016), pp.~621--648.

\bibitem{BurqR}
{\sc N.~Burq and N.~Tzvetkov}, {\em Random data {C}auchy theory for
  supercritical wave equations. {I}. {L}ocal theory}, Invent. Math., 173
  (2008), pp.~449--475.

\bibitem{CarlesJems}
{\sc R.~Carles, E.~Dumas, and C.~Sparber}, {\em Geometric optics and
  instability for {NLS} and {D}avey-{S}tewartson models}, J. Eur. Math. Soc.
  (JEMS), 14 (2012), pp.~1885--1921.

\bibitem{CarlesT}
{\sc R.~Carles and T.~Kappeler}, {\em Norm-inflation with infinite loss of
  regularity for periodic {NLS} equations in negative {S}obolev spaces}, Bull.
  Soc. Math. France, 145 (2017), pp.~623--642.

\bibitem{christ2003ill}
{\sc M.~Christ, J.~Colliander, and T.~Tao}, {\em Ill-posedness for nonlinear
  {S}chr{\"o}dinger and wave equations}.
\newblock arXiv:0311048, 2003.

\bibitem{Elenawave}
{\sc E.~Cordero and F.~Nicola}, {\em Remarks on {F}ourier multipliers and
  applications to the wave equation}, J. Math. Anal. Appl., 353 (2009),
  pp.~583--591.

\bibitem{feichtinger1983modulation}
{\sc H.~G. Feichtinger}, {\em Modulation spaces on locally compact {A}belian
  groups}.
\newblock Technical Report, University of Vienna, 1983, and in “Wavelets and
  Their Applications” (eds. M. Krishna, R. Radha and S. Thangavelu), 99-140,
  Allied Publishers, New Delhi, 2003., 1983.
\newblock Available on
  {researchgate.net}.

\bibitem{JH}
{\sc J.~Forlano and T.~Oh}, {\em Normal form approach to the one-dimensional
  cubic nonlinear {S}chr{\"o}dinger equation in {F}ourier-amalgam spaces,
  \textit{preprint}}.

\bibitem{forlano2019remark}
{\sc J.~Forlano and M.~Okamoto}, {\em
  \href{https://doi.org/10.4310/DPDE.2020.v17.n4.a3}{A remark on norm inflation
  for nonlinear wave equations}}, Dynamics of Partial Differential Equations,
  (2020), pp.~361--381.

\bibitem{grochenig2013foundations}
{\sc K.~Gr\"{o}chenig}, {\em Foundations of time-frequency analysis}, Applied
  and Numerical Harmonic Analysis, Birkh\"{a}user Boston, Inc., Boston, MA,
  2001.

\bibitem{IwabuchiOgawa}
{\sc T.~Iwabuchi and T.~Ogawa}, {\em Ill-posedness for the nonlinear
  {S}chr\"{o}dinger equation with quadratic non-linearity in low dimensions},
  Trans. Amer. Math. Soc., 367 (2015), pp.~2613--2630.

\bibitem{kishimoto2019remark}
{\sc N.~Kishimoto}, {\em \href{http://dx.doi.org/10.3934/cpaa.2019067}{A remark
  on norm inflation for nonlinear {S}chr{\"o}dinger equations}}, Communications
  on Pure \& Applied Analysis, 18 (2019), p.~1375.

\bibitem{Lebeau}
{\sc G.~Lebeau}, {\em Perte de r\'{e}gularit\'{e} pour les \'{e}quations
  d'ondes sur-critiques}, Bull. Soc. Math. France, 133 (2005), pp.~145--157.

\bibitem{LHans}
{\sc H.~Lindblad}, {\em A sharp counterexample to the local existence of
  low-regularity solutions to nonlinear wave equations}, Duke Math. J., 72
  (1993), pp.~503--539.

\bibitem{Lindblad}
{\sc H.~Lindblad and C.~D. Sogge}, {\em On existence and scattering with
  minimal regularity for semilinear wave equations}, J. Funct. Anal., 130
  (1995), pp.~357--426.

\bibitem{OhT}
{\sc T.~Oh}, {\em A remark on norm inflation with general initial data for the
  cubic nonlinear {S}chr\"{o}dinger equations in negative {S}obolev spaces},
  Funkcial. Ekvac., 60 (2017), pp.~259--277.

\bibitem{Reich}
{\sc M.~Reich and M.~Reissig}, {\em Wave equations in modulation spaces---decay
  versus loss of regularity}, in New tools for nonlinear {PDE}s and
  application, Trends Math., Birkh\"{a}user/Springer, Cham, 2019, pp.~371--390.

\bibitem{CV1}
{\sc M.~Ruzhansky, M.~Sugimoto, J.~Toft, and N.~Tomita}, {\em Changes of
  variables in modulation and {W}iener amalgam spaces}, Math. Nachr., 284
  (2011), pp.~2078--2092.

\bibitem{ruzhansky2012modulation}
{\sc M.~Ruzhansky, M.~Sugimoto, and B.~Wang}, {\em Modulation spaces and
  nonlinear evolution equations}, in Evolution equations of hyperbolic and
  {S}chr\"{o}dinger type, vol.~301 of Progr. Math., Birkh\"{a}user/Springer
  Basel AG, Basel, 2012, pp.~267--283.

\bibitem{RuzhanskyTurunen}
{\sc M.~Ruzhansky and V.~Turunen}, {\em Pseudo-differential operators and
  symmetries}, vol.~2 of Pseudo-Differential Operators. Theory and
  Applications, Birkh\"{a}user Verlag, Basel, 2010.
\newblock Background analysis and advanced topics.

\bibitem{TaoW}
{\sc T.~Tao}, {\em Low regularity semi-linear wave equations}, Comm. Partial
  Differential Equations, 24 (1999), pp.~599--629.

\bibitem{wang2007global}
{\sc B.~Wang and H.~Hudzik}, {\em The global {C}auchy problem for the {NLS} and
  {NLKG} with small rough data}, J. Differential Equations, 232 (2007),
  pp.~36--73.

\bibitem{wang2011harmonic}
{\sc B.~Wang, Z.~Huo, C.~Hao, and Z.~Guo}, {\em Harmonic analysis method for
  nonlinear evolution equations. {I}}, World Scientific Publishing Co. Pte.
  Ltd., Hackensack, NJ, 2011.

\bibitem{baoxiang2006isometric}
{\sc B.~Wang, Z.~Lifeng, and G.~Boling}, {\em Isometric decomposition
  operators, function spaces {$E^\lambda_{p,q}$} and applications to nonlinear
  evolution equations}, J. Funct. Anal., 233 (2006), pp.~1--39.

\end{thebibliography}
\end{document}